\documentclass[12pt]{article}
\usepackage{amssymb}

\newtheorem{definition}{Definition}[section]
\newtheorem{theorem}[definition]{Theorem}
\newtheorem{lemma}[definition]{Lemma}
\newtheorem{corollary}[definition]{Corollary}

\newtheorem{note}[definition]{Note}

\def\K{\mathbb K}

\def\Mdf{{\hbox{Mat}}_{d+1}(\K )}

\begin{document}

\title{ \bf Two linear transformations each tridiagonal \\ 
with respect to an eigenbasis of the other; comments on 
the split decomposition
}
\author{Paul Terwilliger}  

\date{}
\maketitle

\centerline{\bf \large Dedicated to Tom Koornwinder on his 60th birthday}
\newpage
\begin{abstract} Let $\K$ denote a field and let
$d$ denote a nonnegative integer. Let
$\mathcal A$ denote a $\K$-algebra isomorphic to 
$\hbox{Mat}_{d+1}(\K)$. An element of $\mathcal A$
is called {\it multiplicity-free} whenever its
eigenvalues are mutually distinct and contained in $\K$.
Let $A$ and $A^*$ denote multiplicity-free
elements in $\mathcal A$. Let $\lbrace E_i\rbrace_{i=0}^d$
(resp. 
$\lbrace E^*_i\rbrace_{i=0}^d$)
denote an ordering of the primitive idempotents of $A$ (resp. $A^*$.)
For $0 \leq i \leq d$ let $\theta_i$ 
(resp. $\theta^*_i$)
denote the eigenvalue
of $A$ (resp. $A^*$)
for $E_i$
(resp. $E^*_i$.) 
Let $V$ denote an
irreducible left $\mathcal A$-module. By a {\it decomposition} 
of $V$ we mean a sequence $\lbrace U_i\rbrace_{i=0}^d$
consisting of 1-dimensional
subspaces of $V$ such that
$V=\sum_{i=0}^d U_i$.
A decomposition $\lbrace U_i\rbrace_{i=0}^d$ of $V$
 is said to be {\it split} 
 (with respect to the orderings
 $\lbrace E_i\rbrace_{i=0}^d$,$\lbrace E^*_i\rbrace_{i=0}^d$) 
whenever both (i)
$(A-\theta_iI)U_i = U_{i+1}$ $(0 \leq i \leq d-1)$,
$(A-\theta_dI)U_d =0$; and
(ii) $(A^*-\theta^*_iI)U_i = U_{i-1}$ $(1 \leq i \leq d)$, 
 $(A^*-\theta^*_0I)U_0 =0$.
We show there exists at most one
decomposition of $V$ which is split
with respect to $\lbrace E_i\rbrace_{i=0}^d$,
$\lbrace E^*_i\rbrace_{i=0}^d$.
We show the following are equivalent:
(i) there exists a decomposition of $V$ which is split
with respect to $\lbrace E_i\rbrace_{i=0}^d,\lbrace E^*_i\rbrace_{i=0}^d$;
(ii)
both
\begin{eqnarray*}
E^*_iAE^*_j = \cases{0, &if $\;i-j> 1$;\cr
\not=0, &if $\;i-j=1$\cr}
 \quad 
&&E_iA^*E_j = \cases{0, &if $\;j-i> 1$;\cr
\not=0, &if $\;j-i=1$\cr}
\end{eqnarray*}
for $0 \leq i,j\leq d$.
We call the sequence
 $(A;A^*;\lbrace E_i\rbrace_{i=0}^d; \lbrace E^*_i\rbrace_{i=0}^d)$ a {\it Leonard system}
whenever both
\begin{eqnarray*}
E^*_iAE^*_j = \cases{0, &if $\;|i-j|> 1$;\cr
\not=0, &if $\;|i-j|=1$\cr}
\quad 
&& E_iA^*E_j = \cases{0, &if $\;|i-j|> 1$;\cr
\not=0, &if $\;|i-j|=1$\cr}
\end{eqnarray*}
for $0 \leq i,j\leq d$.
We show 
 $(A;A^*;\lbrace E_i\rbrace_{i=0}^d; \lbrace E^*_i\rbrace_{i=0}^d)$ is a Leonard system
if and only if both (i)
there exists a decomposition of $V$ which is split with
respect to $\lbrace E_i\rbrace_{i=0}^d,\lbrace E^*_i\rbrace_{i=0}^d$;
(ii) there exists a decomposition of $V$ which is split with
respect to $\lbrace E_{d-i}\rbrace_{i=0}^d,\lbrace E^*_i\rbrace_{i=0}^d$.
We also show 
 $(A;A^*;\lbrace E_i\rbrace_{i=0}^d;\lbrace E^*_i\rbrace_{i=0}^d)$ is a Leonard system
if and only if both (i)
there exists a decomposition of $V$ which is split with
respect to $\lbrace E_i\rbrace_{i=0}^d, \lbrace E^*_i\rbrace_{i=0}^d$;
(ii)
there exists an antiautomorphism
$\dagger$ of $\mathcal A $ such that $A^\dagger=A$ and $A^{*\dagger}=A^*$.
\end{abstract}
\newpage

\section{Leonard pairs and Leonard systems}

\medskip
\noindent
We begin by recalling the notion of a 
 {\it Leonard pair}
\cite{TD00},
\cite{TersubI},
\cite{LS99},
\cite{qSerre},
\cite{LS24},
\cite{conform},
\cite{lsint},
\cite{Terint}.
We will use the following terms.
Let $X$ denote a square matrix.
Then $X$ is called {\it tridiagonal}
whenever each nonzero entry lies on either the diagonal,
the subdiagonal, or the superdiagonal. Assume $X$ is tridiagonal.
Then $X$ is called {\it irreducible} whenever each entry on
the subdiagonal is nonzero and each entry on the superdiagonal is
nonzero.

\medskip
\noindent 
We now define  a  Leonard pair.
For the rest of this paper $\K$ will denote a field.

\begin{definition}
\label{def:lp}
\cite[Definition 1.1]{LS99}
\rm
Let $V$ denote a vector space over $\K$ with finite positive
dimension. 
By a {\it Leonard pair} on $V$, we mean an ordered pair
of linear transformations 
$A: V\rightarrow V$ and
$A^*: V\rightarrow V$ 
 which satisfy both (i), (ii)
below.
\begin{enumerate}
\item There exists a basis for $V$ with
respect to which the matrix representing $A$ is irreducible
tridiagonal and the matrix representing $A^*$ is diagonal.
\item There exists a basis for $V$ with
respect to which the matrix representing $A^*$ is irreducible
tridiagonal and the matrix representing $A$ is diagonal.
\end{enumerate}
\end{definition}

\begin{note}
\rm
According to a common notational convention 
$A^*$ denotes the
conjugate-transpose of $A$. We emphasize we are not using this
convention. In a Leonard pair $A,A^*$ the linear transformations
$A$ and $A^*$ are arbitrary subject to (i), (ii) above.
\end{note}

\noindent Our use of the name ``Leonard pair'' is motivated by
a connection to a theorem of D. Leonard 
\cite[p. 260]{BanIto},
\cite{Leon}
which involves the  $q$-Racah polynomials 
\cite{AWil},
\cite[p. 162]{GR}
and some 
related polynomials of the Askey scheme
\cite{KoeSwa}.
This connection is discussed
in \cite[Appendix A]{LS99} and \cite[Section 16]{LS24}.
See 
\cite{Zhed},
\cite{Grun},
\cite{Koelink3},
\cite{Hjal},
\cite{Zhidd}
for related topics.

\medskip
\noindent 
When working with a Leonard pair, it is often convenient to
consider a closely related and somewhat more abstract concept
called a {\it Leonard system}. In order to define this
we 
recall a few more terms.
Let $d$ denote a nonnegative
integer.
Let 
$\hbox{Mat}_{d+1}(\K)$
denote the $\K$-algebra consisting of all
$d+1$ by $d+1$ matrices which have entries in $\K$. We
index the rows and columns by $0,1,\ldots, d$.
Let $\K^{d+1}$ denote the $\K$-vector space
consisting of all $d+1$ by $1$ matrices which have entries in $\K$.
We index the rows by $0,1,\ldots, d$.
We view $\K^{d+1}$ as a left module for
$\hbox{Mat}_{d+1}(\K)$ under  matrix multiplication.
We observe this module is irreducible.
For the rest of this paper we let $\mathcal A$ 
denote a $\K$-algebra
isomorphic to 
$\hbox{Mat}_{d+1}(\K)$. When we refer to an
$\mathcal A$-module we mean
a left $\mathcal A$-module.
Let $V$
denote an irreducible $\mathcal A$-module.
We remark that $V$ is unique up to isomorphism of 
$\mathcal A$-modules and that $V$ has dimension $d+1$. 
Let $v_0, v_1, \ldots, v_d$ denote a basis for $V$.
For $X
 \in {\mathcal A}$ and for
$Y\in \hbox{Mat}_{d+1}(\K)$,
we say $Y$ 
{\it represents $X$ with respect to $v_0, v_1, \ldots, v_d$}
whenever 
$Xv_j = \sum_{i=0}^d Y_{ij}v_i$ for $0 \leq j\leq d$.
Let $A$ denote an element of $\mathcal A$.
 We say $A $ is {\it
multiplicity-free} whenever it has
$d+1$ distinct  eigenvalues 
in $\K$.
Assume $A$ is multiplicity-free.
Let $\theta_0, \theta_1, \ldots, \theta_d$ denote an ordering of 
the eigenvalues
of $A$, and for $0 \leq i \leq d$  
put
\begin{eqnarray*}
E_i = \prod_{{0 \leq  j \leq d}\atop
{j\not=i}} {{A-\theta_j I}\over {\theta_i-\theta_j}},
\end{eqnarray*}
where $I$ denotes the identity of $\mathcal A$.
We observe
(i) $AE_i = \theta_iE_i$ $(0 \leq i \leq d)$;
(ii) $E_iE_j = \delta_{ij}E_i$ $(0 \leq i,j\leq d)$;
(iii) $\sum_{i=0}^d E_i = I$;
(iv) $A=\sum_{i=0}^d \theta_i E_i$.
Let $\mathcal D$ denote the subalgebra of 
 $\mathcal A$ generated by $A$.
Using
(i)--(iv) 
we find
$E_0, E_1, \ldots, E_d$ is a basis for the 
$\K$-vector space   
$\mathcal D$.
We call $E_i$  the {\it primitive idempotent} of
$A$ associated with $\theta_i$.
It is helpful to think of these primitive idempotents as follows. 
Observe
$$
V = E_0V + E_1V + \cdots + E_dV \qquad \qquad (\hbox{direct sum}).
$$
For $0\leq i \leq d$, $E_iV$ is the (one dimensional) eigenspace of
$A$ in $V$ associated with the 
eigenvalue $\theta_i$, 
and $E_i$ acts  on $V$ as the projection onto this eigenspace. 
We remark that the sequence
$\lbrace A^i |0 \leq i \leq d\rbrace $ is a basis for 
the $\K$-vector space 
$\mathcal D$ 
and that $\prod_{i=0}^d (A-\theta_iI)=0$.
By a {\it Leonard pair in $\mathcal A$}
we mean an ordered pair of elements taken from
$\mathcal A$ which act on $V$
 as a Leonard pair
in the sense of Definition \ref{def:lp}.
We now define a Leonard system.

\begin{definition}
\label{def:ls}
\cite[Definition 1.4]{LS99}
\rm
By a {\it Leonard system} in $\mathcal A$, we mean a sequence
$(A;A^*; $ $\lbrace E_i \rbrace_{i=0}^d; 
\lbrace E^*_i \rbrace_{i=0}^d)$ which satisfies (i)--(v) below.
\begin{enumerate}
\item Each of $A,A^*$ is a multiplicity-free element of $\mathcal A$.
\item $E_0, E_1, \ldots, E_d$ is an ordering of the primitive idempotents
of $A$.
\item $E^*_0, E^*_1, \ldots, E^*_d$ is an ordering of the primitive idempotents
of $A^*$.
\item 
${\displaystyle{
E^*_iAE^*_j = \cases{0, &if $\;|i-j|> 1$;\cr
\not=0, &if $\;|i-j|=1$\cr}
\qquad \qquad 
(0 \leq i,j\leq d).
}}$
\item
${\displaystyle{
 E_iA^*E_j = \cases{0, &if $\;|i-j|> 1$;\cr
\not=0, &if $\;|i-j|=1$\cr}
\qquad \qquad 
(0 \leq i,j\leq d).
}}$
\end{enumerate}
\end{definition}
We
comment on how 
Leonard pairs and Leonard systems are related.
In the following discussion $V$ denotes an irreducible 
 $\mathcal A$-module.
Let 
$(A;A^*;\lbrace E_i \rbrace_{i=0}^d; 
\lbrace E^*_i \rbrace_{i=0}^d)$
denote a Leonard system  in $\mathcal A$.
For $0 \leq i \leq d$ let $v_i$
denote a nonzero vector in $E_iV$.
Then the sequence
$v_0, v_1, \ldots, v_d$
is a basis for $V$ which satisfies 
Definition
\ref{def:lp}(ii). 
For $0 \leq i \leq d$ let $v^*_i$
denote a nonzero vector in $E^*_iV$.
Then the sequence 
$v^*_0, v^*_1, \ldots, v^*_d$
is a basis for $V$ which satisfies 
Definition
\ref{def:lp}(i). By these comments
the pair $A,A^*$ is
a Leonard pair in $\mathcal A$.
Conversely
let $A,A^*$ denote a Leonard pair in $\mathcal A$. 
By \cite[Lemma 1.3]{LS99}
each of $A,A^*$ is multiplicity-free.
 Let $v_0, v_1, \ldots, v_d$ denote a basis for $V$ which satisfies
Definition
\ref{def:lp}(ii). 
For $0 \leq i \leq d$ the vector
$v_i$ is an eigenvector for $A$; 
let $E_i$ denote the corresponding primitive idempotent.
Let $v^*_0, v^*_1, \ldots, v^*_d$ denote a basis for $V$ which satisfies
Definition
\ref{def:lp}(i). 
For $0 \leq i \leq d$ the vector
$v^*_i$ is an eigenvector for $A^*$; 
let $E^*_i$ denote the corresponding 
primitive idempotent.
Then 
$(A;A^*;\lbrace E_i \rbrace_{i=0}^d; 
\lbrace E^*_i \rbrace_{i=0}^d)$ 
is a Leonard system  in $\mathcal A$.
In summary we have the following.

\begin{lemma}
Let $A$ and $A^*$ denote elements in $\mathcal A$. Then the pair $A,A^*$
is a Leonard pair in $\mathcal A$ if and only if
the following (i), (ii) hold.
\begin{enumerate}
\item Each of $A,A^*$ is multiplicity-free.
\item There exists an ordering $E_0, E_1, \ldots, E_d$
of the primitive idempotents of $A$
and there exists an ordering
 $E^*_0, E^*_1, \ldots, E^*_d$
of the primitive idempotents of $A^*$ such that
$(A;A^*;$ $\lbrace E_i \rbrace_{i=0}^d; 
\lbrace E^*_i \rbrace_{i=0}^d)$ is a Leonard system in $\mathcal A$.
\end{enumerate}
\end{lemma}

\noindent Later in this paper we will obtain two characterizations
of a Leonard system. These characterizations are based on
a concept which we call the {\it split decomposition}.
This concept is explained in the next section.

\section{The split decomposition} 

\noindent In  \cite{LS99} we introduced the  split decomposition
for Leonard systems and in \cite{LS24} we discussed this
decomposition in detail.
For our present purposes it is useful to define the split
decomposition in a more general context.
We will refer to the following set-up.

\medskip
\begin{definition} 
\label{def:settup}
\rm
Let $A$ and $A^*$ denote multiplicity-free
elements in $\mathcal A$. Let $E_0, E_1, \ldots, E_d$
denote an ordering of the primitive idempotents of $A$ and
for $0 \leq i \leq d$ let $\theta_i$ denote the eigenvalue
of $A$ for $E_i$.
Let $E^*_0, E^*_1, \ldots, E^*_d$
denote an ordering of the primitive idempotents of $A^*$ and
for $0 \leq i \leq d$ let $\theta^*_i$ denote the eigenvalue
of $A^*$ for $E^*_i$.
We let $\mathcal D$ (resp. ${\mathcal D}^*$) denote the subalgebra
of $\mathcal A$ generated by $A$ (resp. $A^*$.)
We let $V$ denote an irreducible $\mathcal A$-module.
\end{definition}
\noindent 
With reference to Definition
\ref{def:settup}, 
by a {\it decomposition} of $V$ we mean a sequence 
$U_0, U_1, \ldots, U_d$ consisting of 1-dimensional
subspaces of $V$ such that
\begin{eqnarray*}
V=U_0+U_1+\cdots + U_d \qquad \qquad (\hbox{direct sum}).
\end{eqnarray*}
We have a comment.
Let $u_0, u_1, \ldots, u_d$ denote a basis for $V$ 
and for $0 \leq i \leq d$ let $U_i$ denote the subspace of
$V$ spanned by $u_i$.
Then the sequence $U_0, U_1, \ldots, U_d$ is a decomposition
of $V$. Conversely, let $U_0, U_1, \ldots, U_d$ denote a decomposition
of $V$. For $0 \leq i \leq d$ let $u_i$ denote a nonzero vector
in $U_i$. Then $u_0, u_1,\ldots, u_d$ is a basis for $V$.

\begin{definition}
\label{def:sp}
\rm
With reference to Definition
\ref{def:settup}, 
let $U_0, U_1, \ldots, U_d$ denote a decomposition of
$V$. We say this decomposition is {\it split}
 (with respect to the 
orderings $E_0, E_1, \ldots, E_d$ and $E^*_0, E^*_1, \ldots, E^*_d$) 
whenever  both
\begin{eqnarray}
(A-\theta_iI)U_i &=& U_{i+1} \quad (0 \leq i \leq d-1), \qquad  
(A-\theta_dI)U_d =0,
\label{eq:split1}
\\
(A^*-\theta^*_iI)U_i &=& U_{i-1} \quad (1 \leq i \leq d), 
\qquad 
(A^*-\theta^*_0I)U_0 =0.
\label{eq:split2}
\end{eqnarray}
\end{definition}
\noindent
Later in this paper we will obtain two characterizations
of a Leonard system which involve the split decomposition.
For the time being we consider the existence and uniqueness of
the split decomposition. We start with uniqueness.

\begin{lemma}
\label{lem:splitshape}
With reference to
Definition
\ref{def:settup}, the following (i), (ii) hold.
\begin{enumerate}
\item
Assume there exists a decomposition $U_0, U_1, \ldots, U_d$
of $V$ which is split with
respect to 
the orderings
$E_0, E_1, \ldots,$ $ E_d$
and $E^*_0, E^*_1, \ldots, E^*_d$.
Then $U_i=\prod_{h=0}^{i-1}(A-\theta_hI)E^*_0V$ and
$U_i=\prod_{h=i+1}^d(A^*-\theta^*_hI)E_dV$ for $0 \leq i \leq d$. 
\item There exists at most one decomposition
of $V$ which is split with respect to
the orderings $E_0, E_1, \ldots, E_d$
and $E^*_0, E^*_1, \ldots, E^*_d$.
\end{enumerate}
\end{lemma}

\noindent {\it Proof:}
\noindent (i)
From the equation on the right in
(\ref{eq:split2})
we find 
 $U_0=E^*_0V$.
Using this and 
(\ref{eq:split1})
we obtain
$U_i = \prod_{h=0}^{i-1}(A-\theta_hI)E^*_0V $
for $0 \leq i \leq d$.
From the equation on the right in
(\ref{eq:split1})
we find 
 $U_d=E_dV$.
Using this and 
(\ref{eq:split2})
we obtain
$U_i = \prod_{h=i+1}^d(A^*-\theta^*_hI)E_dV $
for $0 \leq i \leq d$.
\\
\noindent (ii) Immediate from (i) above.
\hfill $\Box $
\\

\noindent We 
 turn our attention to
the existence of the split decomposition.
In Section 4 we will give a necessary and sufficient condition for 
this existence.
We will use the following result.

\begin{lemma}
\label{lem:lu}
With reference to 
Definition
\ref{def:settup},
assume there exists a decomposition 
$U_0, U_1, $ $ \ldots, U_d$
of $V$ which is split with
respect to 
the orderings
$E_0, E_1, \ldots,$ $ E_d$
and $E^*_0, E^*_1, \ldots, E^*_d$.
Then the following (i)--(v) hold for  
 $0 \leq i \leq d$.
\begin{enumerate}
\item
$
\sum_{h=0}^i U_h = \sum_{h=0}^i A^hE^*_0V.
$
\item
$
\sum_{h=0}^i U_h = \sum_{h=0}^i E^*_hV.
$
\item
$
\sum_{h=i}^d U_h = \sum_{h=0}^{d-i} A^{*h}E_dV.
$
\item
$
\sum_{h=i}^d U_h = \sum_{h=i}^d E_hV.
$
\item
$U_i = (E^*_0V+E^*_1V+\cdots + E^*_iV)\cap 
 (E_iV+E_{i+1}V+\cdots + E_dV).
$
\end{enumerate}
\end{lemma}
\noindent {\it Proof:}
(i) For $0 \leq j\leq d$
we have
$U_j= \prod_{h=0}^{j-1}(A-\theta_hI)E^*_0V$ 
by Lemma
\ref{lem:splitshape}(i)
so $U_j\subseteq \sum_{h=0}^{j}A^hE^*_0V$.
Apparently
$\sum_{h=0}^i U_h\subseteq \sum_{h=0}^{i}A^hE^*_0V$.
In this inclusion the sum on the left has dimension $i+1$
since 
$U_0, U_1, \ldots, U_d$ is a decomposition.
The sum on the right has dimension at most $i+1$.
Therefore 
$\sum_{h=0}^i U_h= \sum_{h=0}^{i}A^hE^*_0V$.
\\
\noindent (ii)
For $0 \leq j\leq d$
 we have 
$\prod_{h=0}^j(A^*-\theta^*_hI)U_j=0$
by
(\ref{eq:split2}) so
$U_j \subseteq \sum_{h=0}^j E^*_hV$.
Apparently
$\sum_{h=0}^iU_h \subseteq \sum_{h=0}^i E^*_hV$.
In this inclusion
each side has dimension $i+1$ so equality
holds.
\\
\noindent (iii) Similar to the proof of (i) above.  
\\
\noindent (iv) Similar to the proof of (ii) above.  
\\
\noindent (v) Combine (ii), (iv) above.
\hfill $\Box $ \\

\section{Some products}

\noindent 
Our next goal is to display a necessary and sufficient condition
for the existence of the split decomposition.
With reference to Definition
\ref{def:settup}, consider 
the  products
\begin{eqnarray*}
E^*_iAE^*_j, \qquad 
E_iA^*E_j \qquad \qquad (0 \leq i,j\leq d).
\end{eqnarray*}
Our condition has to do with 
which of these products
is 0. In order to motivate our result
we initially consider just one of these products.

\begin{lemma}
\label{lem:Bmeaining}
With reference to Definition 
\ref{def:settup}, 
for $0 \leq  i \leq d$ let $v^*_i$ denote a nonzero vector
in $E^*_iV$ and observe $v^*_0, v^*_1, \ldots, v^*_d$ is
a basis for $V$. Let $B$ denote the matrix in 
$\hbox{Mat}_{d+1}(\K)$
which 
represents
$A$ with respect to this basis, so that
\begin{eqnarray}
Av^*_j = 
 \sum_{i=0}^d B_{ij}v^*_i \qquad \qquad  
(0 \leq j \leq d).
\label{eq:bac}
\end{eqnarray}
Then for  $0 \leq i,j\leq d$ the following are equivalent:
(i) $E^*_iAE^*_j=0$; (ii) 
$B_{ij}=0$.
\end{lemma}
\noindent {\it Proof:}
Let the integers $i,j$ be given.
Observe $E^*_rv^*_s=\delta_{rs}v^*_s$
for $0 \leq r,s\leq d$.
By this and (\ref{eq:bac}) we find
$E^*_iAE^*_jV$ is spanned by $B_{ij}v^*_i$.
The result follows.
\hfill $\Box $ \\

\noindent 
In the next lemma we consider a certain pattern of vanishing
products among the $E^*_iAE^*_j$.
We will use the following notation. Let 
$\lambda$ denote an indeterminate and let 
$\K \lbrack \lambda \rbrack $ denote the $\K$-algebra
consisting of all polynomials in $\lambda$
which have coefficients
in $\K$. Let $f_0, f_1, \ldots, f_d$ denote a sequence
of polynomials taken from $\K\lbrack \lambda \rbrack$.
We say this sequence is {\it graded} whenever 
 $f_i$ has degree exactly $i$ for 
$0 \leq i \leq d$.

\begin{lemma} 
\label{thm:poly}
With reference to Definition \ref{def:settup},
 the following (i)--(iii) are equivalent.
\begin{enumerate}
\item ${\displaystyle{
E^*_iAE^*_j = \cases{0, &if $\;i-j> 1$;\cr
\not=0, &if $\;i-j=1$\cr}
\qquad \qquad 
(0 \leq i,j\leq d)}}$.
\item 
There exists a graded sequence of polynomials
$f_0, f_1, \ldots, f_d$ taken from $\K\lbrack \lambda \rbrack$
such that
$E^*_iV = f_i(A)E^*_0V$ for $0 \leq i\leq d$.

\item For $0 \leq i \leq d$,
\begin{eqnarray}
\sum_{h=0}^iE^*_hV = \sum_{h=0}^iA^hE^*_0V.
\label{eq:long}
\end{eqnarray}
\end{enumerate}
\end{lemma}
\noindent {\it Proof:}
$(i)\Rightarrow (ii)$
For $0 \leq i \leq d$ let $v^*_i$ denote a nonzero vector
in $E^*_iV$ and observe $v^*_0, v^*_1, \ldots, v^*_d$ is
a basis for $V$.
Let $B$ denote the matrix
in 
$\hbox{Mat}_{d+1}(\K)$
which 
represents
$A$ with respect to
this basis.
By  Lemma
\ref{lem:Bmeaining},
\begin{eqnarray}
B_{ij} = \cases{0, &if $\;i-j> 1$;\cr
\not=0, &if $\;i-j=1$\cr}
\qquad \qquad 
(0 \leq i,j\leq d).
\label{eq:bcon}
\end{eqnarray}
Let $f_0, f_1, \ldots, f_d$ denote the polynomials
in $\K\lbrack \lambda \rbrack $
which satisfy $f_0=1$ and 
\begin{eqnarray}
\lambda f_j= \sum_{i=0}^{j+1} B_{ij}f_i \qquad \qquad 
(0 \leq j \leq d-1).
\label{eq:fcon}
\end{eqnarray}
We observe  
$f_i$ has degree exactly
$i$ for $0 \leq i \leq d$  so
the sequence $f_0, f_1, \ldots, f_d$ is graded.
Comparing 
(\ref{eq:bac})  and 
(\ref{eq:fcon}) in light of
(\ref{eq:bcon})
we find $v^*_i=f_i(A)v^*_0$ for
$0 \leq i\leq d$. 
It follows
$E^*_iV = f_i(A)E^*_0V$ for $0 \leq i\leq d$.
\\
\noindent 
$(ii)\Rightarrow (iii)$
For $0 \leq j \leq d$ we have 
$E^*_jV=f_j(A)E^*_0V$. The degree
of $f_j$ is $j$ so
$E^*_jV \subseteq \sum_{h=0}^jA^hE^*_0V$.
Apparently
$\sum_{h=0}^iE^*_hV \subseteq \sum_{h=0}^iA^hE^*_0V$.
In this inclusion the sum on the left
has dimension $i+1$ and
the sum on the right has dimension at most
$i+1$.
Therefore 
$\sum_{h=0}^iE^*_hV = \sum_{h=0}^iA^hE^*_0V$.
\\
\noindent 
$(iii)\Rightarrow (i)$
For $0 \leq i \leq d$ let $V_i$ denote the subspace on
the left or right in 
(\ref{eq:long}).
From the right-hand side of
(\ref{eq:long}) we find
$V_i+AV_i=V_{i+1}$ for $0 \leq i \leq d-1$.
From the left-hand side of 
(\ref{eq:long}) we find
$E^*_rV_s=0$ for $0 \leq s<r\leq d$.
Let $i,j$ denote integers $(0 \leq i,j\leq d)$
and first assume $i-j>1$. We show $E^*_iAE^*_j=0$. 
Observe $E^*_jV\subseteq V_j$ 
and $AV_j\subseteq V_{j+1}$  
so $AE^*_jV\subseteq V_{j+1}$.
However $E^*_iV_{j+1}=0$ since $i-j>1$ 
so $E^*_iAE^*_jV=0$. It follows
 $E^*_iAE^*_j=0$. 
Next we assume $i-j=1$
and show $E^*_iAE^*_j\not=0$.
Suppose 
$E^*_iAE^*_j=0$. 
Then by our previous remarks 
$E^*_iAE^*_h=0$ for $0 \leq h \leq j$. 
By this and since $V_j=\sum_{h=0}^j E^*_hV$ we find 
$E^*_iAV_j=0$.
However
$V_i=V_j+AV_j $ and $E^*_iV_j=0$ so
 $E^*_iV_i=0$. This contradicts the construction 
 so $E^*_iAE^*_j\not=0$.
\hfill $\Box $ 

\begin{corollary}
\label{cor:iso}
With reference to
 Definition \ref{def:settup},
let $v^*_0$ denote a nonzero vector in $E^*_0V$ and
consider the 
$\K$-linear transformation from 
$\mathcal D$ to $V$ which sends
$X$ to $Xv^*_0$ for all $X \in {\mathcal D}$.
Assume the equivalent conditions (i)--(iii) hold
in Lemma
\ref{thm:poly}. Then this linear transformation is an
isomorphism.
\end{corollary}
\noindent {\it Proof:}
Since the $\K$-vector spaces $\mathcal D$ and $V$ have the same dimension
it suffices to show the linear transformation is surjective.
Setting $i=d$ in 
(\ref{eq:long}) we find $V={\mathcal D}v^*_0$. Therefore
the linear transformation
is surjective.
\hfill $\Box $ \\

\noindent Replacing 
 $(A;A^*;\lbrace E_i\rbrace_{i=0}^d; \lbrace E^*_i\rbrace_{i=0}^d)$
by 
 $(A^*;A;\lbrace E^*_{d-i}\rbrace_{i=0}^d; \lbrace E_{d-i}\rbrace_{i=0}^d)$
in Lemma 
\ref{thm:poly} and
Corollary
\ref{cor:iso}
we routinely obtain the following results.

\begin{lemma} 
\label{thm:polydual}
With reference to Definition \ref{def:settup},
the following (i)--(iii) are equivalent.
\begin{enumerate}
\item ${\displaystyle{
E_iA^*E_j = \cases{0, &if $\;j-i> 1$;\cr
\not=0, &if $\;j-i=1$\cr}
\qquad \qquad 
(0 \leq i,j\leq d)}}$.
\item 
There exists a graded sequence of polynomials
$f^*_0, f^*_1, \ldots, f^*_d$ taken from $\K\lbrack \lambda \rbrack$
such that
$E_iV = f^*_{d-i}(A^*)E_dV$ for $0 \leq i\leq d$.
\item  For $0 \leq i \leq d$,
\begin{eqnarray}
\sum_{h=i}^dE_hV = \sum_{h=0}^{d-i}A^{*h}E_dV.
\nonumber
\end{eqnarray}
\end{enumerate}
\end{lemma}

\begin{corollary}
\label{cor:isodual}
With reference to
 Definition \ref{def:settup},
let $v_d$ denote a nonzero vector in $E_dV$ and
consider the 
$\K$-linear transformation from 
${\mathcal D}^*$ to $V$ which sends
$X$ to $Xv_d$ for all $X \in {\mathcal D}^*$.
Assume the equivalent conditions (i)--(iii) hold
in Lemma
\ref{thm:polydual}.
Then this linear transformation is an
isomorphism.
\end{corollary}

\section{The existence of the split decomposition}

\noindent We now display a necessary and sufficient condition for
the existence of the split decomposition.

\begin{theorem}
\label{thm:splitcond}
With reference to 
Definition
\ref{def:settup},
the following (i), (ii) are equivalent.
\begin{enumerate}
\item There exists a decomposition of $V$
which is split 
with respect to the orderings $E_0, E_1, \ldots, $ $ E_d$
and 
 $E^*_0, E^*_1, \ldots, E^*_d$.
\item Both
\begin{eqnarray}
E^*_iAE^*_j &=& \cases{0, &if $\;i-j> 1$;\cr
\not=0, &if $\;i-j=1$\cr}
\qquad \qquad 
(0 \leq i,j\leq d),
\label{eq:cases1}
\\
E_iA^*E_j &=& \cases{0, &if $\;j-i> 1$;\cr
\not=0, &if $\;j-i=1$\cr}
\qquad \qquad 
(0 \leq i,j\leq d).
\label{eq:cases2}
\end{eqnarray}
\end{enumerate}
\end{theorem}
\noindent {\it Proof:}
$(i)\Rightarrow (ii)$ 
By assumption there exists a decomposition
$U_0,U_1, \ldots, U_d$ of $V$ which
 is split with respect to the orderings
 $E_0, E_1, \ldots, E_d$ and $E^*_0, E^*_1, \ldots, E^*_d$. 
For $0 \leq i \leq d$
we have
$
\sum_{h=0}^i U_h = \sum_{h=0}^i A^hE^*_0V
$
by 
 Lemma
\ref{lem:lu}(i)
and 
$
\sum_{h=0}^i U_h = \sum_{h=0}^i E^*_hV 
$
by Lemma
\ref{lem:lu}(ii)
so 
$\sum_{h=0}^i E^*_hV =
\sum_{h=0}^i A^hE^*_0V$.
This gives Lemma \ref{thm:poly}(iii).
Applying that lemma we obtain
 (\ref{eq:cases1}).
For $0 \leq i \leq d$
we have
$
\sum_{h=i}^d U_h = \sum_{h=0}^{d-i} A^{*h}E_dV
$
by 
 Lemma
\ref{lem:lu}(iii)
and 
$
\sum_{h=i}^d U_h = \sum_{h=i}^d E_hV 
$
by Lemma
\ref{lem:lu}(iv)
so 
$\sum_{h=i}^d E_hV =
\sum_{h=0}^{d-i} A^{*h}E_dV$.
This gives Lemma \ref{thm:polydual}(iii).
Applying that lemma we obtain
 (\ref{eq:cases2}).
\\
\noindent
$(ii)\Rightarrow (i)$
For $0 \leq i \leq d$ we define
$\tau_i = \prod_{h=0}^{i-1}(A-\theta_hI)$.
We observe $\tau_0, \tau_1, \ldots, \tau_d$
 is a basis for 
the $\K$-vector space $\mathcal D$.
Let $v^*_0$ denote a nonzero vector in $E^*_0V$.
Observe 
Lemma 
\ref{thm:poly}(i) holds by
(\ref{eq:cases1}) so 
 Corollary
\ref{cor:iso} applies;
by that corollary
$\tau_i v^*_0$ $(0 \leq i \leq d)$  is a basis for
$V$.    
We define $U_i=\hbox{Span}(\tau_iv^*_0)$
for $0 \leq i \leq d$  and observe
$U_0, U_1, \ldots, U_d$ is a decomposition
of $V$. We show this decomposition is split with 
respect to 
$E_0, E_1, \ldots, E_d$ and 
$E^*_0, E^*_1, \ldots, E^*_d$.
To do this we 
show the sequence $U_0, U_1, \ldots, U_d$ 
satisfies 
(\ref{eq:split1}) and 
(\ref{eq:split2}).
Concerning (\ref{eq:split1}),
 from the construction
$(A-\theta_iI)\tau_i = \tau_{i+1}$
for $0 \leq i \leq d-1$
and  $(A-\theta_dI)\tau_d=0$.
Applying both sides of these equations to $v^*_0$ we find
$(A-\theta_iI)U_i = U_{i+1}$
for $0 \leq i \leq d-1$
and  $(A-\theta_dI)U_d=0$.
We have now shown 
(\ref{eq:split1}).
Concerning (\ref{eq:split2}),
this will follow if we can show
(a) $(A^*-\theta^*_iI)U_i \subseteq \sum_{h=0}^{i-1} U_h$ 
for $0 \leq i \leq d$; 
 (b) $(A^*-\theta^*_iI)U_i \subseteq \sum_{h=i-1}^d U_h$
for $1 \leq i \leq d$;
(c) 
 $(A^*-\theta^*_iI)U_i\not=0$
for $1 \leq i \leq d$.
We begin with (a).
For $0 \leq j \leq d$ the 
elements $\lbrace \tau_h|0 \leq h \leq j\rbrace$ and the
elements 
$\lbrace A^h|0 \leq h \leq j\rbrace$ span the 
same 
subspace of $\mathcal D$.
Therefore $\sum_{h=0}^j U_h = \sum_{h=0}^j A^hE^*_0V$.
We mentioned
Lemma
\ref{thm:poly}(i) holds
so Lemma \ref{thm:poly}(iii) holds;
therefore 
$ \sum_{h=0}^j E^*_hV=
 \sum_{h=0}^j A^hE^*_0V$
so  
$\sum_{h=0}^jU_h =\sum_{h=0}^j E^*_hV$.
Observe $(A^*-\theta^*_iI)
\sum_{h=0}^i E^*_hV 
=\sum_{h=0}^{i-1} E^*_hV 
$ for $0 \leq i \leq d$.
Combining these comments we find
 $(A^*-\theta^*_iI)U_i \subseteq \sum_{h=0}^{i-1} U_h$ 
for $0 \leq i \leq d$.  
We now have
 (a).
Next we prove (b).
From the construction,
for $0 \leq j\leq d$ we have
$\prod_{h=j}^d (A-\theta_hI)\tau_j=0$
so 
$\prod_{h=j}^d (A-\theta_hI)U_j=0$.
From this we find
$
U_j\subseteq \sum_{h=j}^d E_hV$.
Apparently
$
\sum_{h=i}^dU_h\subseteq \sum_{h=i}^d E_hV$ for
$0 \leq i\leq d$.
By this and since 
 $U_0, U_1, \ldots, U_d$ 
is a decomposition we find
$\sum_{h=i}^dU_h =\sum_{h=i}^d E_hV$ for
$0 \leq i\leq d$.
From 
(\ref{eq:cases2}) 
we find $A^*E_jV \subseteq \sum_{h=j-1}^d E_hV$ for
$1 \leq j\leq d$.
Therefore 
$(A^*-\theta^*_jI)\sum_{h=j}^d E_hV \subseteq 
\sum_{h=j-1}^d E_hV$ for $1 \leq j\leq d$.
From these comments we find
$(A^*-\theta^*_jI)U_j \subseteq 
\sum_{h=j-1}^dU_h$  for $1 \leq j\leq d$.
We now have (b).
Next we show (c). Suppose there exists
an integer $i$ $(1 \leq i\leq d)$ such that
$(A^*-\theta^*_iI)U_i=0$.
We assume $i$ is maximal subject to this.
We obtain a contradiction as follows.
For $i < j\leq d$ we find
$(A^*-\theta^*_jI)U_j\subseteq U_{j-1} $
by (a), (b).
In this inclusion
the left-hand side is nonzero and the right-hand side has
dimension 1 so we have equality.
We mentioned earlier $(A-\theta_dI)U_d=0$ so
$U_d=E_dV$.
Apparently 
 $U_j=\prod_{h=j+1}^d (A^*-\theta^*_hI)E_dV$ for
$i\leq j\leq d$. 
In particular
 $U_i=\prod_{h=i+1}^d (A^*-\theta^*_hI)E_dV$.
Combining this with
$(A^*-\theta^*_iI)U_i=0$ we obtain
$0=\prod_{h=i}^d (A^*-\theta^*_hI)E_dV$.
Let $v_d$ denote a nonzero vector in $E_dV$
and observe 
$0=\prod_{h=i}^d (A^*-\theta^*_hI)v_d$.
This is inconsistent with 
 Corollary
\ref{cor:isodual}
and the fact that
$0\not=\prod_{h=i}^d (A^*-\theta^*_hI)$.
We now have a contradiction and
(c) is proved.
Combining (a)--(c)
we obtain
(\ref{eq:split2}).
We have shown the decomposition $U_0,U_1, \ldots, U_d$
satisfies 
(\ref{eq:split1}), 
(\ref{eq:split2}). Applying Definition
\ref{def:sp} we find 
$U_0,U_1, \ldots, U_d$
is split with respect to 
the orderings
$E_0, E_1, \ldots, E_d$
and 
 $E^*_0, E^*_1, \ldots, E^*_d$.
\hfill $\Box $ \\

\section{Two characterizations of a Leonard system}

In this section we obtain two characterizations of a
Leonard system, both of which involve the 
split decomposition.
We will first state the characterizations,
then prove a few lemmas, and then prove 
the characterizations.
Our first characterization
is stated as follows.

\begin{theorem}
\label{thm:char1}
With reference to Definition
\ref{def:settup},
the sequence
$(A;A^*;\lbrace E_i\rbrace_{i=0}^d; \lbrace E^*_i\rbrace_{i=0}^d)$
 is a Leonard system if
and only if both  (i), (ii) hold below.
\begin{enumerate}
\item There exists a decomposition of $V$ which is split with
respect to the orderings
$E_0, E_1, \ldots, $ $ E_d$ and 
$E^*_0, E^*_1, \ldots, E^*_d$.
\item There exists a decomposition of $V$ which is split with
respect to the orderings
$E_d, E_{d-1}, $ $ \ldots, E_0$ and 
$E^*_0, E^*_1, \ldots, E^*_d$.
\end{enumerate}
\end{theorem}

\noindent In order to state our second characterization we 
recall a
definition.
Let $\sigma :
{\mathcal A}\rightarrow {\mathcal A}$
denote any map. We call $\sigma$ 
an
{\it antiautomorphism}
of $\mathcal A$ 
whenever $\sigma$ is an
isomorphism of $\K$-vector spaces  
and  $(XY)^\sigma = Y^\sigma X^\sigma $ for
all $X, Y \in {\mathcal A}$.
For example assume 
${\mathcal A}=\hbox{Mat}_{d+1}(\K)$.
 Then $\sigma$ is an antiautomorphism
of $\mathcal A$ if and only if there exists an invertible
$R \in {\mathcal A}$ such that
$X^\sigma = R^{-1}X^tR$ for all $X \in {\mathcal A}$,
where $t$ denotes transpose. This follows from the
Skolem-Noether Theorem
\cite[Cor. 9.122]{rotman}.

\medskip
\noindent We now state our second characterization of
a Leonard system.

\begin{theorem}
\label{thm:char2}
With reference to Definition
\ref{def:settup},
the sequence 
 $(A;A^*;\lbrace E_i\rbrace_{i=0}^d; \lbrace E^*_i\rbrace_{i=0}^d)$
 is a Leonard system if
and only if both (i), (ii) hold below.
\begin{enumerate}
\item There exists a decomposition of $V$ which is split with
respect to the orderings
$E_0, E_1, \ldots, $ $ E_d$ and 
$E^*_0, E^*_1, \ldots, E^*_d$.
\item There exists an antiautomorphism
$\dagger$ of $\mathcal A $ such that $A^\dagger=A$ and $A^{*\dagger}=A^*$.
\end{enumerate}
\end{theorem}

\noindent We now prove some lemmas which we will use to
obtain
Theorem \ref{thm:char1} and 
Theorem \ref{thm:char2}.
We have a preliminary remark.
With reference to 
Definition
\ref{def:settup}, we consider the following four conditions:
\begin{eqnarray}
&&\quad E^*_iAE^*_j = \cases{0, &if $\; i-j > 1$;\cr
\not=0, &if $\; i-j = 1$\cr}
\qquad \qquad 
(0 \leq i,j\leq d),
\label{eq:fourth3}
\\
&&\quad E^*_iAE^*_j = \cases{0, &if $\; j-i > 1$;\cr
\not=0, &if $\; j-i = 1$\cr}
\qquad \qquad 
(0 \leq i,j\leq d),
\label{eq:fourth4}
\\
&&\quad E_iA^*E_j = \cases{0, &if $\; i-j > 1$;\cr
\not=0, &if $\; i-j = 1$\cr}
\qquad \qquad 
(0 \leq i,j\leq d),
\label{eq:fourth1}
\\
&&\quad E_iA^*E_j = \cases{0, &if $\; j-i > 1$;\cr
\not=0, &if $\; j-i = 1$\cr}
\qquad \qquad 
(0 \leq i,j\leq d).
\label{eq:fourth2}
\end{eqnarray}
We observe 
 $(A;A^*;\lbrace E_i\rbrace_{i=0}^d; \lbrace E^*_i\rbrace_{i=0}^d)$
is a Leonard system if and only if
each of 
(\ref{eq:fourth3})--(\ref{eq:fourth2}) holds.

\begin{lemma}
With reference to Definition
\ref{def:settup},
 assume
conditions (\ref{eq:fourth3}) and
(\ref{eq:fourth4}) hold.
Then $A,E^*_0$ together generate $\mathcal A$.
Moreover $A,A^*$ together generate $\mathcal A$.
\end{lemma}
\noindent {\it Proof:}
Examining the proof of
\cite[Lemma 3.1]{LS24} we find that the 
elements $A^rE^*_0A^s$ $(0\leq r,s\leq d)$ form a basis
for the $\K$-vector space $\mathcal A$.
It follows that $A,E^*_0$ together generate $\mathcal A$.
The elements $A,A^*$ together generate $\mathcal A$ since
$E^*_0$ is a polynomial
in $A^*$.
\hfill $\Box $ \\

\begin{lemma}
\label{lem:dagexist}
With reference to Definition
\ref{def:settup},
 assume conditions  
(\ref{eq:fourth3})  and 
(\ref{eq:fourth4}) hold.
Then 
there exists a unique antiautomorphism $\dagger$ of $\mathcal A$ such
that $A^\dagger = A$ and $A^{*\dagger}=A^*$.
Moreover $X^{\dagger \dagger}=X$ for all $X \in {\mathcal A}$. 
\end{lemma}
\noindent {\it Proof:}
Concerning the existence of $\dagger$,
for $0 \leq i \leq d$ let $v^*_i$ denote a nonzero element
of $E^*_iV$
and recall $v^*_0, v^*_1, \ldots, v^*_d$ is
a basis for $V$. For $X \in {\mathcal A}$ let
$X^\flat$ denote the matrix in 
$\hbox{Mat}_{d+1}(\K)$ which represents 
$X$ with respect to the basis 
$v^*_0, v^*_1, \ldots, v^*_d$.
We observe $\flat : {\mathcal A}\rightarrow 
\hbox{Mat}_{d+1}(\K)$ 
is an isomorphism of $\K$-algebras.
We abbreviate $B=A^\flat$ and 
$B^*=A^{*\flat}$.
We observe $B$ is irreducible tridiagonal and
$B^*=\hbox{diag}(\theta^*_0, \theta^*_1, \ldots, \theta^*_d)$.
Let $D$ denote the diagonal matrix 
in $\Mdf$ which has $ii$ entry
\begin{eqnarray*}
D_{ii}= \frac
{B_{01}B_{12}\cdots B_{i-1,i}}
{B_{10}B_{21}\cdots B_{i,i-1}}
\qquad \qquad (0 \leq i \leq d).
\end{eqnarray*}
It is routine to verify $D^{-1}B^tD=B$.
Each of $D, B^*$ is diagonal 
so
$DB^* = 
B^*D$;
 also $B^{*t}=B^*$ so
$D^{-1}B^{*t}D=B^*$. 
Let $\sigma :\Mdf \rightarrow \Mdf$ denote the 
map which satisfies
$X^\sigma = D^{-1}X^tD$  for all
$X \in \Mdf$. We observe
$\sigma$ is an antiautomorphism
of $\Mdf$ such that
$B^\sigma=B$ and $B^{*\sigma}=B^*$.
We define the map 
 $\dagger :{\mathcal A} \rightarrow {\mathcal A}$ 
to be the composition
$\dagger:=\flat \sigma \flat^{-1}$.
We observe
$\dagger $ is an antiautomorphism
of $\mathcal A$ such that $A^\dagger = A$ and 
$A^{*\dagger}=A^*$.
We have now shown there exists an antiautomorphism
$\dagger $
of $\mathcal A$ such that
 $A^\dagger=A$ and $A^{*\dagger}=A^*$.
This antiautomorphism is
unique since $A, A^*$ together generate $\mathcal A$.
The map $X\rightarrow X^{\dagger \dagger}$
is an isomorphism of $\K$-algebras from 
$\mathcal A$ to itself. This map is the identity
since 
$A^{\dagger \dagger}=A$, 
$A^{*\dagger \dagger}=A^*$,
and since $A,A^*$ together generate $\mathcal A$.
\hfill $\Box $ \\

\begin{lemma}
\label{lem:eistab}
With reference to Definition
\ref{def:settup},
assume there exists an
antiautomorphism $\dagger$ of $\mathcal A$ such
that $A^\dagger = A$ and $A^{*\dagger}=A^*$.
Then $E_i^\dagger = E_i$ and 
$E_i^{*\dagger} = E^*_i$ 
for $0 \leq i \leq d$.
\end{lemma}
\noindent {\it Proof:}
 Recall $E_i$ (resp. $E^*_i$) is a polynomial
in $A$ (resp. $A^*$) for $0 \leq i \leq d$.
\hfill $\Box $ \\

\begin{lemma}
\label{lem:flip}
With reference to Definition
\ref{def:settup},
assume there exists an
antiautomorphism $\dagger$ of $\mathcal A$ such
that $A^\dagger = A$ and $A^{*\dagger}=A^*$.
Then for $0 \leq i,j\leq d$,
(i)
$E^*_iAE^*_j=0$ if and only if
$E^*_jAE^*_i=0$; and 
(ii)
$E_iA^*E_j=0$ if and only if
$E_jA^*E_i=0$.
\end{lemma}
\noindent {\it Proof:}
By Lemma
\ref{lem:eistab} and since $\dagger $ is an antiautomorphism,
\begin{eqnarray*}
(E^*_iAE^*_j)^\dagger &=& E^*_jAE^*_i \qquad \qquad (0 \leq i,j\leq d).
\end{eqnarray*}
Assertion (i) follows since $\dagger :{\mathcal A}\rightarrow {\mathcal A} $ is 
a bijection.
To obtain (ii) interchange the roles of $A$ and 
$A^*$ in the proof of (i).
\hfill $\Box $ \\

\begin{lemma}
\label{lem:3gives4}
With reference to Definition
\ref{def:settup},
assume at least three of 
(\ref{eq:fourth3})--(\ref{eq:fourth2})
hold. Then each of
(\ref{eq:fourth3})--(\ref{eq:fourth2})
hold; in other words
 $(A;A^*;\lbrace E_i\rbrace_{i=0}^d; \lbrace E^*_i\rbrace_{i=0}^d)$
  is a Leonard system.
\end{lemma}
\noindent {\it Proof:}
Interchanging $A$ and $A^*$ if necessary,
we may assume without loss of generality that 
(\ref{eq:fourth3}) and 
(\ref{eq:fourth4}) hold.
By Lemma \ref{lem:dagexist}
 there exists an 
antiautomorphism $\dagger $ of $\mathcal A$ such 
that $A^\dagger=A$ and $A^{*\dagger}=A^*$.
By assumption at least one of 
(\ref{eq:fourth1}),
(\ref{eq:fourth2}) holds.
Combining  this 
with Lemma
\ref{lem:flip} we find 
(\ref{eq:fourth1}),
(\ref{eq:fourth2}) both hold.
The result follows.
\hfill $\Box $ \\

\noindent 
We are now ready to prove Theorem
\ref{thm:char1}.
\medskip

\noindent {\it Proof of 
 Theorem
\ref{thm:char1}:}
By Theorem
\ref{thm:splitcond} we find
(i) holds if and only if each of
(\ref{eq:fourth3}),
(\ref{eq:fourth2}) holds.
Applying Theorem
\ref{thm:splitcond} again, this time with
 $(A;A^*;\lbrace E_i\rbrace_{i=0}^d; \lbrace E^*_i\rbrace_{i=0}^d)$
replaced by
 $(A;A^*;\lbrace E_{d-i}\rbrace_{i=0}^d; \lbrace E^*_i\rbrace_{i=0}^d)$,
we find
(ii) holds if and only if each of
(\ref{eq:fourth3}),
(\ref{eq:fourth1}) holds.
Suppose 
 $(A;A^*;\lbrace E_i\rbrace_{i=0}^d; \lbrace E^*_i\rbrace_{i=0}^d)$
is a Leonard system.
Then each of 
(\ref{eq:fourth3})--(\ref{eq:fourth2}) holds.
In particular
each of 
(\ref{eq:fourth3}),
(\ref{eq:fourth1}),
(\ref{eq:fourth2}) holds so
(i), (ii) hold by our above remarks.
Conversely suppose (i), (ii) hold.
Then each of  
(\ref{eq:fourth3}),
(\ref{eq:fourth1}), 
(\ref{eq:fourth2}) 
holds.
At least three of
(\ref{eq:fourth3})--(\ref{eq:fourth2}) hold
so 
 $(A;A^*;\lbrace E_i\rbrace_{i=0}^d; \lbrace E^*_i\rbrace_{i=0}^d)$
is a Leonard system
 by Lemma
\ref{lem:3gives4}.
\hfill $\Box $ \\

\noindent
We are now ready to prove Theorem
\ref{thm:char2}.

\medskip
\noindent {\it Proof of Theorem 
\ref{thm:char2}:}
First assume
 $(A;A^*;\lbrace E_i\rbrace_{i=0}^d; \lbrace E^*_i\rbrace_{i=0}^d)$
is a Leonard system.
Then (i) holds by 
Theorem \ref{thm:char1} and 
 (ii) holds by Lemma
\ref{lem:dagexist}.
Conversely assume (i), (ii) hold.
Combining (i) and 
Theorem
\ref{thm:splitcond} we 
obtain 
(\ref{eq:fourth3}),
(\ref{eq:fourth2}).
Combining this with (ii) and using
Lemma \ref{lem:flip}
we obtain 
(\ref{eq:fourth4}),
(\ref{eq:fourth1}).
Now each of 
(\ref{eq:fourth3})--(\ref{eq:fourth2}) holds
so 
 $(A;A^*;\lbrace E_i\rbrace_{i=0}^d; \lbrace E^*_i\rbrace_{i=0}^d)$
is a Leonard system.
\hfill $\Box $ \\

\noindent We would like to emphasize the following fact.

\begin{theorem}
\label{thm:lpdag}
Let $A,A^*$ denote a Leonard pair in $\mathcal A$. Then there
exists a unique antiautomorphism $\dagger $ 
of $\mathcal A$ such that $A^\dagger=A$ and $A^{*\dagger}=A^*$.
Moreover $X^{\dagger \dagger}=X$ for all $X \in {\mathcal A}$.
\end{theorem}
\noindent {\it Proof:}
 Since $A,A^*$ is a Leonard pair
there exists an ordering $E_0, E_1, \ldots, E_d$ of the primitive
idempotents of $A$ and an ordering 
 $E^*_0, E^*_1, \ldots, E^*_d$ of the primitive
idempotents of $A^*$ such that 
 $(A;A^*;\lbrace E_i\rbrace_{i=0}^d; \lbrace E^*_i\rbrace_{i=0}^d)$
 is a Leonard system.
These orderings satisfy 
(\ref{eq:fourth3})--(\ref{eq:fourth2}).
In particular
(\ref{eq:fourth3}), (\ref{eq:fourth4}) are satisfied
so the result follows by Lemma 
\ref{lem:dagexist}.
\hfill $\Box $ \\

\noindent We finish this section with a comment.

\begin{lemma}
\label{lem:lp}
With reference to 
Definition
\ref{def:settup}, 
 assume there exists a decomposition
of $V$
which is split with respect to
the orderings
$E_0, E_1, \ldots, E_d$ and
$E^*_0, E^*_1, \ldots, E^*_d$.
Then the following 
(i), (ii) are equivalent.
\begin{enumerate}
\item 
The pair
$A,A^*$ is a Leonard
pair.
\item 
The sequence
 $(A;A^*;\lbrace E_i\rbrace_{i=0}^d; \lbrace E^*_i\rbrace_{i=0}^d)$
is a Leonard system.
\end{enumerate}
\end{lemma}
\noindent {\it Proof:}
$(i)\Rightarrow (ii)$ 
We assume
 there exists a decomposition
of $V$
which is split with respect to
the orderings
$E_0, E_1, \ldots, E_d$ and
$E^*_0, E^*_1, \ldots, E^*_d$.
Therefore each of  
(\ref{eq:fourth3}),
 (\ref{eq:fourth2}) holds by
 Theorem
\ref{thm:splitcond}.
Since $A,A^*$ is a Leonard pair there exists
an antiautomorphism $\dagger$ of $\mathcal A$ such that
$A^\dagger=A$ and $A^{*\dagger}=A^*$.
Applying Lemma
\ref{lem:flip} we find
each of (\ref{eq:fourth4}),
 (\ref{eq:fourth1}) holds. Now each of
(\ref{eq:fourth3})--(\ref{eq:fourth2}) 
holds so 
 $(A;A^*;\lbrace E_i\rbrace_{i=0}^d; \lbrace E^*_i\rbrace_{i=0}^d)$
is a Leonard system.
\\
$(ii)\Rightarrow (i)$ Clear.
\hfill $\Box $ \\

\section{The two characterizations in terms of matrices}

\noindent In this section we restate
Theorem \ref{thm:char1} and 
Theorem \ref{thm:char2}
in terms of matrices.
We first set some notation.
With reference to 
Definition
\ref{def:settup},
suppose there exists a  decomposition 
$U_0, U_1, \ldots, U_d$
of $V$ which is split with respect to
the orderings
$E_0, E_1, \ldots, E_d$
and 
 $E^*_0, E^*_1, \ldots, E^*_d$.
Pick an integer $i$ $(1 \leq i \leq d)$. By
(\ref{eq:split2}) we find
$(A^*-\theta^*_iI)U_i= U_{i-1}$ and by
(\ref{eq:split1}) we find
$(A-\theta_{i-1}I)U_{i-1}= U_i$.
Apparently $U_i$ is an eigenspace for  
$(A-\theta_{i-1}I)
(A^*-\theta^*_iI)$ and the corresponding eigenvalue
is a nonzero element of $\K$.
Let us denote this eigenvalue by $\varphi_i$.
We call $\varphi_1, \varphi_2, \ldots, \varphi_d$
the {\it split sequence} for $A,A^*$ with respect to
the orderings
$E_0, E_1, \ldots, E_d$
and 
 $E^*_0, E^*_1, \ldots, E^*_d$.
The split sequence has the following interpretation.
For $0 \leq i \leq d$ let $u_i$ denote a nonzero vector
in $U_i$ and recall $u_0,u_1, \ldots, u_d$ is a basis for 
$V$. We normalize the $u_i$ so that
$(A-\theta_iI)u_i = u_{i+1}$ for $0 \leq i \leq d-1$.
With respect to the basis $u_0, u_1,\ldots, u_d$
the matrices which represent $A$
and $A^*$ are as follows.
\begin{eqnarray*}
A:\; \left(
\begin{array}{c c c c c c}
\theta_0 & & & & & {\mathbf 0} \\
1 & \theta_1 &  & & & \\
& 1 & \theta_2 &  & & \\
& & \cdot & \cdot &  &  \\
& & & \cdot & \cdot &  \\
{\mathbf 0}& & & & 1 & \theta_d
\end{array}
\right),
\qquad  A^*: 
\left(
\begin{array}{c c c c c c}
\theta^*_0 &\varphi_1 & & & & {\mathbf 0} \\
 & \theta^*_1 & \varphi_2 & & & \\
&  & \theta^*_2 & \cdot & & \\
& &  & \cdot & \cdot &  \\
& & &  & \cdot & \varphi_d \\
{\mathbf 0}& & & &  & \theta^*_d
\end{array}
\right).
\end{eqnarray*}
Motivated by this 
we consider the following set-up.

\begin{definition} 
\label{def:fset}
Let $d$ denote a nonnegative integer. 
Let $A$ and $A^*$ denote matrices in $\Mdf$ of the form
\begin{eqnarray*}
A=\left(
\begin{array}{c c c c c c}
\theta_0 & & & & & {\mathbf 0} \\
1 & \theta_1 &  & & & \\
& 1 & \theta_2 &  & & \\
& & \cdot & \cdot &  &  \\
& & & \cdot & \cdot &  \\
{\mathbf 0}& & & & 1 & \theta_d
\end{array}
\right),
\qquad  A^*= 
\left(
\begin{array}{c c c c c c}
\theta^*_0 &\varphi_1 & & & & {\mathbf 0} \\
 & \theta^*_1 & \varphi_2 & & & \\
&  & \theta^*_2 & \cdot & & \\
& &  & \cdot & \cdot &  \\
& & &  & \cdot & \varphi_d \\
{\mathbf 0}& & & &  & \theta^*_d
\end{array}
\right),
\end{eqnarray*}
where
\begin{eqnarray*}
&&\theta_i\not=\theta_j, \qquad 
\theta^*_i\not=\theta^*_j \quad \hbox{if} \quad i\not=j
\qquad   \qquad (0\leq i,j\leq d),
\\
&& \qquad \varphi_i \not=0 \qquad \qquad (1 \leq i\leq d).
\end{eqnarray*}
We observe $A$ (resp. $A^*$) is multiplicity-free,
with eigenvalues $\theta_0, \theta_1,\ldots, \theta_d$
(resp. $\theta^*_0, \theta^*_1,\ldots,$ $ \theta^*_d$.)
For $0 \leq i \leq d$ we let $E_i$  (resp. $E^*_i$) denote
the primitive idempotent for $A$ (resp. $A^*$) associated
with $\theta_i$ (resp. $\theta^*_i$.)
\end{definition}

\noindent We have some comments.
With reference to
Definition
\ref{def:fset},
for $0 \leq i \leq d$ let $u_i$ denote
the vector in 
$\K^{d+1}$ which has ith entry 1 and all other entries 0.
We observe $u_0, u_1, \ldots, u_d$ is a basis for
$\K^{d+1}$.
From the form of $A$ we have
$(A-\theta_iI)u_i=u_{i+1}$ for $0 \leq i \leq d-1$
and 
 $(A-\theta_dI)u_d=0$.
From the form of $A^*$ we have
$(A^*-\theta^*_iI)u_i=\varphi_iu_{i-1}$ for $1 \leq i \leq d$
and 
 $(A^*-\theta^*_0I)u_0=0$.
For $0 \leq i \leq d$ let $U_i$ denote the subspace of
$\K^{d+1}$ spanned by $u_i$. Then $U_0, U_1, \ldots, U_d$
is a decomposition of 
$\K^{d+1}$.
This decomposition satisfies
$(A-\theta_iI)U_i=U_{i+1}$ for $0 \leq i \leq d-1$
and 
 $(A-\theta_dI)U_d=0$.
Similarly
$(A^*-\theta^*_iI)U_i=U_{i-1}$ for $1 \leq i \leq d$
and $(A^*-\theta^*_0I)U_0=0$.
In other words the decomposition
$U_0, U_1, \ldots, U_d$ is split with respect to
the orderings
$E_0, E_1, \ldots, E_d$ and
$E^*_0, E^*_1, \ldots, E^*_d$. We observe $\varphi_1, \varphi_2,
\ldots, \varphi_d$ is the corresponding split sequence for
$A,A^*$.
We now consider when is the pair $A,A^*$ a Leonard pair.
We begin with a remark. 

\begin{lemma}
\label{lem:mention}
With reference to
Definition
\ref{def:fset}, the following (i), (ii) are equivalent. 
\begin{enumerate}
\item 
The pair $A,A^*$ is a Leonard pair.
\item 
The sequence  
 $(A;A^*;\lbrace E_i\rbrace_{i=0}^d; \lbrace E^*_i\rbrace_{i=0}^d)$
 is a Leonard system.
\end{enumerate}
\end{lemma}
\noindent {\it Proof:}
We mentioned  there exists a 
 decomposition of $\K^{d+1}$ 
which is split with respect to
the orderings
$E_0, E_1, \ldots, E_d$ and
$E^*_0, E^*_1, \ldots, E^*_d$.
Therefore 
Lemma \ref{lem:lp} applies and
the result follows.
\hfill $\Box $ \\

\noindent We now give a matrix version of Theorem
\ref{thm:char1}.

\begin{theorem}
\label{thm:bbcc}
Referring to Definition
\ref{def:fset},
the following (i), (ii) are equivalent.
\begin{enumerate}
\item The pair
$A,A^*$ is a Leonard pair.
\item There exists an invertible 
$G \in \hbox{Mat}_{d+1}(\K)$ and there exists
 nonzero $\phi_i
\in \K$ $(1 \leq i \leq d)$ 
 such that
\begin{eqnarray*}
G^{-1}AG = 
\left(
\begin{array}{c c c c c c}
\theta_d & & & & & {\mathbf 0} \\
1 & \theta_{d-1} &  & & & \\
& 1 & \theta_{d-2} &  & & \\
& & \cdot & \cdot &  &  \\
& & & \cdot & \cdot &  \\
{\mathbf 0}& & & & 1 & \theta_0
\end{array}
\right),
\quad 
G^{-1}A^*G = 
\left(
\begin{array}{c c c c c c}
\theta^*_0 &\phi_1 & & & & {\mathbf 0} \\
 & \theta^*_1 & \phi_2 & & & \\
&  & \theta^*_2 & \cdot & & \\
& &  & \cdot & \cdot &  \\
& & &  & \cdot & \phi_d \\
{\mathbf 0}& & & &  & \theta^*_d
\end{array}
\right).
\end{eqnarray*}
\end{enumerate}
Suppose (i), (ii) hold. Then the sequence
$\phi_1, \phi_2, \ldots, \phi_d$
is the split sequence for $A,A^*$ associated with
the orderings 
$E_d, E_{d-1}, \ldots, E_0$ and
$E^*_0, E^*_1, \ldots, E^*_d$.
\end{theorem}
\noindent {\it Proof:}
$(i)\Rightarrow (ii)$
The sequence  
 $(A;A^*;\lbrace E_i\rbrace_{i=0}^d; \lbrace E^*_i\rbrace_{i=0}^d)$
is a Leonard system
by 
Lemma 
\ref{lem:mention}.
By Theorem \ref{thm:char1} there
 exists a decomposition
of $\K^{d+1}$ which is split
with respect to 
the orderings
$E_d, E_{d-1}, $ $\ldots, E_0$ and $E^*_0, E^*_1, \ldots, E^*_d$.
Let 
$V_0, V_1, \ldots, V_d$
denote this decomposition.
By the definition of a split decomposition
we have
$(A-\theta_{d-i}I)V_i = V_{i+1}$ for
$0 \leq i \leq d-1$ and
$(A-\theta_0I)V_d =0$.
Moreover
$(A^*-\theta^*_iI)V_i = V_{i-1}$ for $1 \leq i \leq d$ 
and 
$(A^*-\theta^*_0I)V_0 =0$.
For $0 \leq i \leq d$ let $v_i$ denote a nonzero vector
in $V_i$ and observe $v_0, v_1, \ldots, v_d$ is a basis
for $\K^{d+1}$.
We normalize the $v_i$ so that $(A-\theta_{d-i}I)v_i = v_{i+1}$
for $0 \leq i \leq d-1$.
Let $\phi_1, \phi_2, \ldots, \phi_d$ denote the 
split sequence for $A,A^*$ with respect to 
the orderings
$E_d, E_{d-1}, \ldots, E_0$
and $E^*_0, E^*_1, \ldots, E^*_d$.
Then $\phi_i\not=0$ $(1 \leq i \leq d)$ and moreover
 $(A^*-\theta^*_iI)v_i =
\phi_i v_{i-1}$ $(1 \leq i \leq d)$,
$(A^*-\theta^*_0I)v_0 =0$.
Let $G$ denote the matrix in 
$\hbox{Mat}_{d+1}(\K)$ which has column $i$ equal to $v_i$
for $0 \leq i \leq d$. We observe $G$ is invertible.
Moreover the matrices
 $G^{-1}AG$ and $G^{-1}A^*G$ have  the form shown above.
\\
\noindent 
$(ii)\Rightarrow (i)$
We show
 $(A;A^*;\lbrace E_i\rbrace_{i=0}^d; \lbrace E^*_i\rbrace_{i=0}^d)$
is a Leonard system.
In order to do this we apply Theorem
\ref{thm:char1}.
In the paragraph after  Definition 
\ref{def:fset} we mentioned
there exists a 
 decomposition of $\K^{d+1}$ 
which is split with respect to
the orderings
$E_0, E_1, \ldots, E_d$ and
$E^*_0, E^*_1, \ldots, E^*_d$.
Therefore 
Theorem \ref{thm:char1}(i)
holds.
We show 
 Theorem
\ref{thm:char1}(ii) holds.
For $0 \leq i \leq d$ let 
$v_i$ denote column $i$ of $G$ and observe $v_0, v_1, \ldots, v_d$
is a basis for $\K^{d+1}$. From the form of $G^{-1}AG$ we find
$(A-\theta_{d-i}I)v_i = v_{i+1}$ for $0 \leq i\leq d-1$ and
$(A-\theta_{0}I)v_d =0$.
From the form of $G^{-1}A^*G$ we find
$(A^*-\theta^*_iI)v_i = \phi_iv_{i-1}$ for $1 \leq i\leq d$ and
$(A^*-\theta^*_{0}I)v_0 =0$.
For $0 \leq i \leq d$ let 
$V_i$ denote the subspace of $\K^{d+1}$ 
spanned
by $v_i$. Then
$V_0, V_1, \ldots, V_d$ is a decomposition
of 
$\K^{d+1}$.
Also 
$(A-\theta_{d-i}I)V_i = V_{i+1}$ for $0 \leq i\leq d-1$ and
$(A-\theta_{0}I)V_d =0$.
Moreover
$(A^*-\theta^*_iI)V_i = V_{i-1}$ for $1 \leq i\leq d$ and
$(A^*-\theta^*_{0}I)V_0 =0$.
Apparently 
$V_0, V_1, \ldots, V_d$ is 
split with respect to 
the orderings
$E_d, E_{d-1}, \ldots, E_0$
and $E^*_0, E^*_1, \ldots, E^*_d$.
Now 
 Theorem \ref{thm:char1}(ii) holds; applying that theorem
we find
 $(A;A^*;\lbrace E_i\rbrace_{i=0}^d; \lbrace E^*_i\rbrace_{i=0}^d)$
is a Leonard system. In particular $A,A^*$ is a Leonard pair.
\\
\noindent Assume
(i), (ii) both hold. From the proof of
$(ii)\Rightarrow (i)$ we find
that for $1 \leq i \leq d$,
$\phi_i$ is the eigenvalue of
$(A-\theta_{d-i+1}I)(A^*-\theta^*_iI)$ associated with $V_i$.
Therefore  $\phi_1, \phi_2, \ldots, \phi_d$
is the split sequence for $A,A^*$ associated with
the orderings
$E_d, E_{d-1}, \ldots, E_0$ and 
$E^*_0, E^*_1, \ldots, E^*_d$.
\hfill $\Box $ \\

\noindent We now give a matrix 
version of Theorem
\ref{thm:char2}.

\begin{theorem}
\label{thm:char2m}
Referring to Definition
\ref{def:fset},
the following (i), (ii) are equivalent.
\begin{enumerate}
\item The pair
$A,A^*$ is a Leonard pair.
\item
There exists an invertible  $H \in \Mdf$ such that
\begin{eqnarray*}
H^{-1}A^tH = A, \qquad \qquad H^{-1}A^{*t}H=A^*.
\end{eqnarray*}
\end{enumerate}
\end{theorem}
\noindent {\it Proof:}
$(i)\Rightarrow (ii)$
By Theorem 
\ref{thm:lpdag}
there exists an antiautomorphism $\dagger$ of
$\hbox{Mat}_{d+1}(\K)$ such that $A^\dagger=A$ and
$A^{*\dagger}=A^*$.
Since $\dagger$ is an antiautomorphism
there exists an invertible $H \in
\hbox{Mat}_{d+1}(\K)$ 
such that $X^\dagger=H^{-1}X^tH$ for all 
$X
\in \hbox{Mat}_{d+1}(\K)$. 
Setting $X=A$ we have
$H^{-1}A^tH = A$. Setting $X=A^*$ we have
$H^{-1}A^{*t}H=A^*$.
\\
\noindent
$(ii)\Rightarrow (i)$
We show 
 $(A;A^*;\lbrace E_i\rbrace_{i=0}^d; \lbrace E^*_i\rbrace_{i=0}^d)$
 is a Leonard system.
In order to do this we apply Theorem
\ref{thm:char2}.
In the paragraph after  Definition 
\ref{def:fset} we mentioned
there exists a 
 decomposition of $\K^{d+1}$ 
which is split with respect to
the orderings
$E_0, E_1, \ldots, E_d$ and
$E^*_0, E^*_1, \ldots, E^*_d$.
Therefore 
Theorem \ref{thm:char2}(i)
holds.
Let $\dagger
:\hbox{Mat}_{d+1}(\K)\rightarrow
\hbox{Mat}_{d+1}(\K)
$
 denote the map 
 which satisfies
$X^\dagger=H^{-1}X^tH$ for all $X \in 
\hbox{Mat}_{d+1}(\K)$.
Then $\dagger$ is an antiautomorphism of
$\hbox{Mat}_{d+1}(\K)$ such that
$A^\dagger=A$ and $A^{*\dagger}=A^*$.
Now
Theorem \ref{thm:char2}(ii) holds;
 applying that theorem
we find 
 $(A;A^*;\lbrace E_i\rbrace_{i=0}^d; \lbrace E^*_i\rbrace_{i=0}^d)$
 is a Leonard system.
In particular $A,A^*$ is a Leonard pair.
\hfill $\Box $ \\

\section{Remarks}

Referring to Definition
\ref{def:fset},
presumably 
condition (ii) of 
Theorem \ref{thm:bbcc} or 
Theorem \ref{thm:char2m}
can be translated into a 
 condition on the entries of
$A$ and $A^*$.
We obtained such a
condition in 
\cite{LS99};
we cite it here for the sake
of completeness.

\begin{theorem}
\cite[Corollary 14.2]{LS99}
With reference to Definition
\ref{def:fset},
the pair $A,A^*$ is a Leonard pair if and only if
there exists nonzero $\phi_i \in \K$ $(1 \leq i \leq d)$
 such that (i)--(iii) hold below.
\begin{enumerate}
\item $\varphi_i = \phi_1 \sum_{h=0}^{i-1} 
\frac{\theta_h-\theta_{d-h}}{\theta_0-\theta_d}
+ (\theta^*_i-\theta^*_0)(\theta_{i-1}-\theta_d)
\qquad \qquad (1 \leq i \leq d)$.
\item $\phi_i = \varphi_1 \sum_{h=0}^{i-1} 
\frac{\theta_h-\theta_{d-h}}{\theta_0-\theta_d}
+ (\theta^*_i-\theta^*_0)(\theta_{d-i+1}-\theta_0)
\qquad \qquad (1 \leq i \leq d)$.
\item The expressions
\begin{eqnarray*} 
\frac{\theta_{i-2}-\theta_{i+1}}{\theta_{i-1}-\theta_i},
\qquad 
\frac{\theta^*_{i-2}-\theta^*_{i+1}}{\theta^*_{i-1}-\theta^*_i}
\end{eqnarray*}
are equal and independent of $i$ for $2 \leq i \leq d-1$.
\end{enumerate}
Suppose (i)--(iii) hold. Then $\phi_1, \phi_2,\ldots, \phi_d$
is the split sequence for $A,A^*$ with respect to 
the orderings
$E_d, E_{d-1}, \ldots, E_0$ and 
$E^*_0, E^*_1, \ldots, E^*_d$.
\end{theorem}

\section{Acknowledgement}
The author would like to thank Brian Curtin,
Eric Egge, 
Mark MacLean,
Arlene Pascasio,
and Chih-wen Weng
for giving the manuscript a close
reading and offering many valuable suggestions.

\noindent Paul Terwilliger \hfil\break
Department of Mathematics \hfil\break
University of Wisconsin \hfil\break
480 Lincoln Drive \hfil\break
Madison, Wisconsin, 53706 USA 
\hfil\break
email: terwilli@math.wisc.edu \hfil\break


\begin{thebibliography}{10}

\bibitem{AWil}
R.~Askey and J.A.~Wilson.
\newblock A set of orthogonal polynomials that 
generalize the {R}acah coefficients or $6-j$ symbols. 
\newblock {\em SIAM J. Math. Anal.}, 10:1008--1016, 1979. 

\bibitem{BanIto}
E.~Bannai and T.~Ito.
\newblock {\em Algebraic Combinatorics {I}: Association Schemes}.
\newblock Benjamin/Cummings, London, 1984.

\bibitem{GR}
G.~Gasper and M.~Rahman.
\newblock {\em Basic hypergeometric series}.
\newblock Encyclopedia of Mathematics and its Applications, 35.
\newblock Cambridge University Press, Cambridge, 1990.


\bibitem{Zhed}
Ya.~Granovskii, I.~Lutzenko, and A. Zhedanov.   
\newblock  Mutual integrability, quadratic algebras, and dynamical symmetry.
\newblock  {\em Ann. Physics}, 217(1):1--20, 1992.


\bibitem{Grun}
F.~A.~Grunbaum and L.~Haine.
\newblock A $q$-version of a theorem of Bochner.
\newblock {\em J. Comput. Appl. Math.}, 68(1-2):103--114, 1996.


\bibitem{TD00}
T.~Ito, K.~Tanabe, and P.~Terwilliger.
\newblock Some algebra related to ${P}$- and ${Q}$-polynomial association
  schemes. In
\newblock {\em Codes and Association Schemes (Piscataway NJ, 1999)}, 
DIMACS Ser. Discrete Math. Theoret. Comput. Sci., 56:167--192,
Amer.
Math. Soc., Providence RI, 2001.


\bibitem{KoeSwa}
R.~Koekoek and R.~Swarttouw.
\newblock {\em The Askey-scheme of hypergeometric orthogonal polyomials and its
  $q$-analog}, volume 98-17 of {\em Reports of the faculty of Technical
  Mathematics and Informatics}.
\newblock Delft, The Netherlands, 1998.
Available at
 \newblock{ \tt http://aw.twi.tudelft.nl/{\~{}}koekoek/research.html}



\bibitem{Koelink3}
H.~T. Koelink.
\newblock Askey-{W}ilson polynomials and the quantum ${\rm {s}{u}}(2)$ group:
  survey and applications.
\newblock {\em Acta Appl. Math.}, 44(3):295--352, 1996.

\bibitem{Leon}
D.~Leonard.
\newblock Orthogonal polynomials, duality, and association
schemes.
\newblock {\em SIAM J. Math. Anal.}, 13(4):656--663, 1982.


\bibitem{Hjal}
H.~Rosengren.
\newblock {\em Multivariable orthogonal polynomials as coupling
coefficients for Lie and quantum algebra representations.}
\newblock Centre for Mathematical Sciences, Lund University, Sweden,
1999.

\bibitem{rotman}
J.~J.~Rotman.
\newblock{\em Advanced modern algebra.}
\newblock{Prentice Hall, Saddle River NJ 2002}.

\bibitem{TersubI}
P.~Terwilliger.
\newblock The subconstituent algebra of an association scheme I. 
\newblock {\em J. Algebraic Combin.}, 1(4):363--388, 1992.

\bibitem{LS99}
P.~Terwilliger.
\newblock Two linear transformations each tridiagonal with respect to an
  eigenbasis of the other.
\newblock {\em Linear Algebra Appl.},  330:149--203, 2001. 

\bibitem{qSerre}
P.~Terwilliger.
\newblock Two relations that generalize the q-Serre relations and the 
Dolan-Grady relations. In
\newblock {\em
Physics and Combinatorics 1999 (Nagoya)}, 377--398, 
World Sci. Publishing, River Edge, NJ, 2001. 


\bibitem{LS24}
P.~Terwilliger.
\newblock  Leonard pairs from 24 points of view.
\newblock {\em Rocky Mountain J. Math.}, 32(2):827--888, 2002. 

\bibitem{conform}
P.~Terwilliger.
\newblock Two linear transformations each tridiagonal with respect to an
  eigenbasis of the other: the $TD$-$D$ and the $LB$-$UB$ canonical form.
Preprint.

\bibitem{lsint}
P.~Terwilliger.
\newblock Introduction to Leonard pairs.
\newblock {\em J. Comput. Appl. Math.} 
(OPSFA Rome 2001), 153(2):463--475, 2003.

\bibitem{Terint}
P.~Terwilliger.
\newblock An Introduction to {L}eonard pairs and 
	     {L}eonard systems.
  \newblock {\em S\=urikaisekikenky\=usho K\=oky\=uroku},
  \newblock 
   (1109):67--79, 1999.
    \newblock Algebraic combinatorics  (Kyoto, 1999).


\bibitem{Zhidd}
A.~S. Zhedanov.
\newblock ``{H}idden symmetry'' of {A}skey-{W}ilson polynomials.
\newblock {\em Teoret. Mat. Fiz.}, 89(2):190--204, 1991.


\end{thebibliography}
\end{document}